\documentclass[14pt]{article}

\usepackage{natbib}
\usepackage{amsmath}
\usepackage{amssymb,times}
\usepackage{amscd}
\usepackage[margin=2cm]{geometry}
\usepackage{helvet}
\usepackage{stmaryrd}
\usepackage{siunitx}
\usepackage{bigstrut}
\usepackage{graphicx}
\usepackage{varwidth}
\usepackage{capt-of}
\usepackage{color}
\usepackage{subcaption}

\usepackage[]{sfmath}
\newcommand{\dd}[2]{\frac{\diff#1}{\diff #2}} 
\def\MM#1{\boldsymbol{#1}}
\newcommand{\pp}[2]{\frac{\partial #1}{\partial #2}} 
\DeclareMathOperator{\diff}{d\!}
\DeclareMathOperator{\Tr}{Tr}
\DeclareMathOperator{\Hdiv}{H(div)}

\newcommand{\correction}[1]{{#1}}

\title{Higher-order compatible finite element schemes for the nonlinear rotating
  shallow water equations on the sphere} \author{J. Shipton, T. H. Gibson and C. J. Cotter}

\begin{document}

\maketitle

\begin{abstract}
  We describe a compatible finite element discretisation for the
  shallow water equations on the rotating sphere, concentrating on
  integrating consistent upwind stabilisation into the
  framework. Although the prognostic variables are velocity and layer
  depth, the discretisation has a diagnostic potential vorticity that
  satisfies a stable upwinded advection equation through a
  Taylor-Galerkin scheme; this provides a mechanism for dissipating
  enstrophy at the gridscale whilst retaining optimal order
  consistency. We also use upwind discontinuous Galerkin schemes for
  the transport of layer depth. These transport schemes are
  incorporated into a semi-implicit formulation that is facilitated by
  a hybridisation method for solving the resulting mixed Helmholtz
  equation.  {We demonstrate that our discretisation achieves the
  expected second order convergence and provide results from some
  standard rotating sphere test problems.}
\end{abstract}

\section{Introduction}

The development of new numerical discretisations based on finite
element methods is being driven by the need for more flexibility in
mesh geometry. The scalability bottleneck arising from the
latitude-longitude grid means that weather and climate model
developers are searching for numerical discretisations that are stable
and accurate on pseudo-uniform grids without sacrificing properties of
conservation, balance and wave propagation that are important for
accurate atmosphere modelling on the scales relevant to weather and
climate \citep{staniforth2012horizontal}.  There is also ongoing
interest in adaptively refined meshes as a way of seamlessly coupling
global scale and local scale atmosphere simulations, as well as
dynamic adaptivity or even moving meshes; using these meshes requires
numerical methods that can remain stable and accurate on multiscale
meshes. Further, there is an interest in using higher-order
spaces to try to offset the inhomogeneity in the error due to using
grids that break rotational symmetry.

Compatible finite element methods are a form of mixed finite element
methods (meaning that different finite element spaces are used for
different fields) that allow the exact representation of the standard
vector calculus identities div-curl=0 and curl-grad=0.  This
necessitates the use of H(div) finite element spaces for
velocity, such as Raviart-Thomas and Brezzi-Douglas-Marini, and
discontinuous finite element spaces for pressure (stable pairing of
velocity and pressure space relies on the existence of bounded
commuting projections from continuous to discrete spaces, as detailed
in \citet{boffi2013mixed}, for example). The main reason for choosing
compatible finite element spaces is that they have a discrete
Helmholtz decomposition of the velocity space; this means that there
is a clean separation between divergence-free and rotational velocity
fields.  \citet{cotter2012mixed} used this decomposition to
demonstrate that compatible finite element discretisations for the
linear shallow water equations on arbitrary grids satisfy the basic
conservation, balance and wave propagation properties listed in
\citet{staniforth2012horizontal}. In particular, it was shown that the
discretisation has a geostrophic balancing pressure for every velocity field in the divergence-free subspace of the H(div) finite element space. A survey
of the stability and approximation properties of compatible finite
element spaces is provided in \citet{natale2016compat}, including a
proof of the absence of spurious inertial oscillations.

The challenge of building atmosphere models using compatible finite
elements is that there is no freedom to select finite element spaces
in order to ensure good representation of the nonlinear equations
(such as conservation, or accurate advection, for example), because
the choice has already been made to satisfy linear requirements. In
the case of the rotating shallow water equations, the use of
discontinuous finite element spaces for the layer depth field
encourages us to use upwind discontinuous Galerkin methods, to solve
the continuity equation describing layer depth transport.

The nonlinearity in the momentum/velocity equation is more
challenging. In \citet{mcrae2014energy}, the energy-enstrophy
conserving formulation of \citet{arakawa1981potential} was extended to
compatible finite element methods. This extension is closely related
to C-grid methods for the shallow water equations on more general
meshes in \citet{ringler2010unified,thuburn2012framework}. Following
these approaches, the compatible finite element formulation, which has
velocity and height as prognostic variables, has a diagnostic
potential vorticity that satisfies a conservation equation that is
implied by the prognostic dynamics for velocity and height. A finite
element exterior calculus structure in this formulation was exposed in
\citet{cotter2014finite}, which also provided an alternative
formulation based around low-order finite element methods on dual
grids. In \citet{thuburn2015primal}, the close relationship of the
dual grid formulation to finite volume methods was exploited to obtain
a stable discretisation of the nonlinear shallow water equations on
the sphere where the finite element formulation of the wave dynamics
was coupled with high-order finite volume methods for the layer depth
and prognostic potential vorticity fields. The essential idea is to
select a particular stable accurate finite volume scheme for the
diagnostic potential vorticity, and to then find the update for the
prognostic velocity which implies it. In this paper we address the
issue of extending this idea to higher-order finite element spaces,
for which there is no analogue of the dual grid spaces. This means
that we must return to the formulation of \citet{mcrae2014energy},
where the potential vorticity is stored in a continuous finite element
space. We then seek stable accurate higher-order discretisations
of the potential vorticity equation using continuous finite element
methods that make it possible to find the corresponding update for
prognostic velocity. It turns out that this is indeed possible for
advection methods from the SUPG/Taylor-Galerkin family of methods.

Finally, we show how these discretisations can be embedded within a
semi-implicit time-integration scheme. We again follow the
formulation in \citet{thuburn2015primal}, in which advection terms are
obtained from explicit time integration methods applied using the
(iterative) velocity at time level $n+1/2$. The linear system solved
during each nonlinear iteration for the corrections to the field
values also requires attention. The standard approach of eliminating
velocity to solve a Helmholtz problem for the correction to the layer
depth is problematic because the inverse velocity mass matrix is
dense. We instead use a hybridised formulation where one solves for
the Lagrange multipliers that enforce normal continuity of the
velocity field \citep[for example]{boffi2013mixed}.

In section \ref{sec: SW model} we describe the shallow water model,
including the spatial and temporal discretisation; we present finite
element spaces that satisfy the properties outlined above and provide
details of how to construct such spaces on the sphere and describe
advection schemes for both discontinuous and continuous fields as
required. In section \ref{sec: results} we present the results of
applying our scheme to some of the standard set of test cases for
simulation of the rotating shallow water equations on the sphere as
described in \citet{williamson1992standard} and
\citet{galewsky2004initial}. Section \ref{sec: summary} provides a
summary and brief outlook.

\section{{The shallow water model}}
\label{sec: SW model}

\subsection{{Shallow water equations}}

We begin with the vector invariant form of the nonlinear shallow water
equations on a two dimensional surface $\Omega$ embedded in three
dimensions,
\begin{eqnarray}
\MM{u}_t + (\zeta +f)\MM{u}^\perp + \nabla \left(g(D+b)+\frac{1}{2}|\MM{u}|^2\right) &=& 0, \label{eq:momentum}\\
D_t + \nabla\cdot(\MM{u}D) &=& 0, \label{eq:mass}
\end{eqnarray}
where $\MM{u}$ is the horizontal velocity, $D$ is the layer depth, $b$
is the height of the lower boundary, $g$ is the gravitational
acceleration, $f$ is the Coriolis parameter and
$\zeta=\nabla^\perp\cdot\MM{u}:=(\MM{k}\times\nabla)\cdot\MM{u}$ is
the vorticity, $\MM{u}^\perp = \MM{k}\times\MM{u}$, $\MM{k}$ is the
normal to the surface $\Omega$, and where the $\nabla$ and
$\nabla\cdot$ operators are defined intrinsically on the
surface. These equations have the important property that the shallow
water potential vorticity (PV)
\begin{equation}
q=\frac{\zeta +f}{D}
\end{equation}
satisfies a local conservation law,
\begin{equation}
\pp{}{t}(Dq) + \nabla \cdot(\MM{u}qD) = 0.
\end{equation}
This can be seen by applying $\nabla^\perp\cdot$ to equation
\eqref{eq:momentum}. Equation \eqref{eq:mass} then implies that
$q$ is constant along characteristics moving with the flow
velocity $\MM{u}$, \emph{i.e.},
\begin{equation}
\label{eq:pv advection}
\pp{q}{t} + \left(\MM{u}\cdot\nabla\right)q=0.
\end{equation}
Numerical discretisations that preserve some aspects of these
properties have been demonstrated to be very successful at obtaining
long time integrations of the rotating shallow water equations on
sphere in the quasi-geostrophic flow regime. This is partially because
they provide a way to avoid discretising the vector advection term
$(\MM{u}\cdot\nabla)\MM{u}$ in the velocity equation directly; instead
one can choose a suitable stable, accurate and conservative scalar
advection scheme for the potential vorticity (treated as a diagnostic
variable with $\MM{u}$ and $D$ being prognostic variables) and use it
to diagnose a form of the $(\zeta+f)\MM{u}^\perp$ term in equation
\eqref{eq:momentum} that leads to stable advection of $\MM{u}$. These
ideas were introduced in the compatible finite element context in
\citet{mcrae2014energy} in order to obtain energy-enstrophy conserving
discretisations; here we concentrate on stable, accurate and
conservative advection of $q$, improving on the low-order APVM
stabilisation suggested there. We also replace the centred
discretisation of equation \eqref{eq:mass} with stable and accurate
Discontinuous Galerkin (DG) advection schemes for $D$, and show how
these can be incorporated into the PV conserving formulation.

\subsection{{Spatial discretisation}}
\subsubsection{Finite element spaces}
\label{FESpaces}

In this section we shall summarise the properties we require from our
finite element spaces and the operators between them. We start with
the space $H(\textrm{div})$ of square integrable velocity fields,
whose divergence in also square integrable. The condition that the
discrete velocity belongs to the finite element subspace
$\mathbb{V}_1\subset H(\textrm{div})$ means that the velocity must
have continuous normal components across element edges. Having chosen
$\mathbb{V}_1$, we select a finite element space $\mathbb{V}_2\subset
L^2$
such that 
\begin{equation}
\{\nabla\cdot\MM{w} : \MM{w}\in\mathbb{V}_1\} \subset \mathbb{V}_2.
\end{equation}
This necessarily requires that $\mathbb{V}_2$ is a discontinuous
space.
We also define a space $\mathbb{V}_0\subset H^1$ consisting of
continuous fields $\gamma$ such that
$\MM{k}\times\nabla\psi\in\mathbb{V}_1$, where the curl
$\MM{k}\times\nabla$, henceforth written as $\nabla^\perp$, maps from
$\mathbb{V}_0$ onto the kernel of $\nabla\cdot$ in $\mathbb{V}_1$.

The proof that the mixed finite element discretisation of the linear
shallow water equations has steady geostrophic modes relies on the
existence of a discrete Helmholtz decomposition for the velocity field
\citep{cotter2012mixed}. As described in \citet{arnold2006finite},
this decomposition exists if the following diagram commutes with
bounded projections $\pi_1$, $\pi_2$, $\pi_3$,
\begin{equation}\begin{CD}
\label{eq:commutative}
  H^1 @> \nabla^\perp >> H(\textrm{div}) @> \nabla\cdot
 >> L^2 \\
  @VV{\pi_0}V @VV{\pi_1}V @VV{\pi_2}V \\
\mathbb{V}_0 @>\nabla^\perp >> \mathbb{V}_1 @> \nabla\cdot >> \mathbb{V}_2
\end{CD}\end{equation}

\noindent that is, the result of applying an operator to the continuous field
and projecting into the discrete space is the same as the result of
first projecting the field into the discrete space and then applying
the operator.

\citet{cotter2012mixed} reviewed several sets of finite element spaces
that satisfy these requirements, together with further requirements on
the degree-of-freedom (DOF) ratios between $\mathbb{V}_1$ and
$\mathbb{V}_2$ that are necessary to exclude the possibility of
spurious mode branches in the dispersion relation for the linear
shallow water equations. In this paper with shall present results from
choosing $\mathbb{V}_1$ to be the first order Brezzi-Douglas-Marini
($\textrm{BDM}_2$) function space on triangles which requires that
$\mathbb{V}_0$ is $\textrm{P}_3$, the space of piecewise continuous
cubic functions. $\mathbb{V}_2$ is the space $\textrm{P}_1^{DG}$
of piecewise linear functions that can be discontinuous at the element
boundaries.


\subsubsection{Constructing the finite element spaces on the sphere}

In order to implement the finite element method, we need to expand our
variables in terms of suitable basis functions and compute integrals
of combinations of these basis functions over elements. This is done
by defining a reference element on which these integrals can be
calculated and mappings from each reference element  to the
physical element, which we describe in this section. A more
detailed and general exposition, plus description of how these
mappings are implemented in the FEniCS project, is provided in
\cite{rognes2013automating}. In the rest of this section, hatted
quantities refer to those defined on the reference element.

We start by defining finite element spaces $\mathbb{V}_i(\hat{e})$,
$i=0,1,2$, on the reference element $\hat{e}$. These are constructed
from polynomials in the usual manner for the chosen spaces as
described in \cite{boffi2013mixed}, for example. Then, for each
element $e$, we construct a polynomial mapping
$\MM{g}_e:\hat{e}\mapsto \mathbb{R}^3$, such that
\begin{equation}
\hat{\Omega} = \cup_{i=1}^{N_e} \MM{g}_{e_i}(\hat{e}), \quad
\MM{g}_{e_1}(\hat{e}) \cap \MM{g}_{e_2}(\hat{e}) = \emptyset, \, 
\mbox{ if }e_1\neq e_2,
\end{equation}
where $\hat{\Omega}$ is the computational domain, which is a piecewise
polynomial approximation to the sphere domain $\Omega$. Then, we use
$\MM{g}_e$ to relate functions in our finite element spaces on
$\hat{\Omega}$, restricted to each element $e$ in $\hat{\Omega}$, to
functions in $\mathbb{V}_i(\hat{e})$.

For $\mathbb{V}_0$, $\mathbb{V}_2$, these functions are simply related
by function composition, \emph{i.e.},
\begin{equation}
\psi \in \mathbb{V}_0 \implies \psi|_e\circ\MM{g}_e \in \mathbb{V}_0(\hat{e}),
\quad 
\phi \in \mathbb{V}_2 \implies \phi|_e\circ\MM{g}_e \in \mathbb{V}_2(\hat{e}),
\label{eq:psi-phi}
\end{equation}

For $\mathbb{V}_1$, we use the Piola transformation {\citep{boffi2013mixed}}

\begin{equation}
\MM{u}|_e\left(\MM{g}_e(\hat{\MM{x}})\right)=\frac{J_e\hat{\MM{u}}}{\det J_e},
\end{equation}
where $\MM{u}|_e$ is the restriction of $\MM{u}$ to element
$e$, and
\begin{equation}
J_e=\pp{\MM{g}_e}{\hat{\MM{x}}}
\end{equation}
is the Jacobian of $\MM{g}$ {in element $e$}, to ensure that $\MM{u}$ is
tangent to $\Omega$. The Piola transformation has the crucial property
that

\begin{equation}
\nabla\cdot\MM{u}
\left(\MM{g}_e(\hat{\MM{x}})\right) = 
\frac{\hat{\nabla}\cdot\hat{\MM{u}}(\hat{\MM{x}})}
{\det J_e\left(\hat{\MM{x}}\right)}.
\end{equation}
When the discrete space is made up of flat elements, as in
\citet{cotter2012mixed} and \citet{mcrae2014energy}, $\det J_e$ is
constant. This guarantees that if
$\hat{\nabla}\cdot\hat{\MM{u}}\in\hat{\mathbb{V}}_1$ then
$\nabla\cdot\MM{u}\in\mathbb{V}_1$. However, for meshes made up of
general quadrilaterals (i.e.\ cubed sphere meshes such as those in
\citet{putman2007finite}) or high-order, curved triangles, the mapping
$\MM{g}$ is no longer affine and the determinant of its Jacobian is no
longer constant. This means that in general
$\nabla\cdot\MM{u}\notin\mathbb{V}_1$. This clearly violates the
commutative diagram property \eqref{eq:commutative}. There are 3 ways
to remedy the situation. 
\begin{enumerate}
\item Modify the mapping for $\mathbb{V}_2$ to become
\begin{equation}
\phi \in \mathbb{V}_2 \implies \det J_e \phi|_e\circ\MM{g}_e \in \mathbb{V}_2(\hat{e}).
\end{equation}
This option is the choice that is most consistent with the {finite element exterior calculus} methodology.
\item Modify the mapping for $\mathbb{V}_1$ so that the factor of
  $\det J_e$ in $\nabla\cdot\MM{u}|_e\circ \MM{g}_e$ is replaced by
  the element average of $\det J_e$. This approach was described in
  \citet{boffi2009some} for the case of lowest order Raviart-Thomas
  elements; the extension to other H(div) elements is straightforward
  but the construction is quite complicated.
\ \item Replace the divergence operator $\nabla\cdot$ in the commutative
  diagram by $\tilde{\nabla}\cdot$ to be the $L_2$ projection $\pi_2$ of
$\nabla\cdot$ into $\mathbb{V}_1$, defined by
\begin{equation}
\label{newdiv}
\int_\Omega\phi\tilde{\nabla}\cdot\MM{u} \diff V = \int_\Omega\phi\nabla\cdot\MM{u} \diff V, \quad \forall \phi \in \mathbb{V}_2.
\end{equation}
By defining a new divergence operator in this way we immediately
recover the required commutation property, since $\pi_2$ appears in
the commutative diagram \eqref{eq:commutative} and is a projection.
This is a generalisation of the ``rehabilitation'' technique described
for lowest order Raviart-Thomas elements applied to mixed elliptic
problems in \cite{bochev2008rehabilitation}. In fact, as Bochev and
Ridzal noticed, introduction of the $\tilde{\nabla}\cdot$ operator
does not require any changes to a code implementation for mixed
elliptic operators since the $\tilde{\nabla}\cdot$ operator only
appears in an inner product with a test function from $\mathbb{V}_2$,
and hence can be safely replaced by $\nabla\cdot$ there. We shall 
see that this continues to be the case in our nonlinear shallow water
formulation.
\end{enumerate}
In general, Option 1 is the most mathematically elegant
choice. However, since our software implementation is based upon the
Firedrake finite element library \citep{Rathgeber2016, Luporini2016}
which is already used to develop DG schemes that use the
transformations in \eqref{eq:psi-phi}, Option 3 was much less
pervasive through the code base. We will investigate the differences
between Options 1 and 3 in future work.

\citet{Arnold14} showed that there is a potential problem with loss of
consistency with compatible finite elements on quadrilateral and cubic
curvilinear cells; it is likely that this problem also can be
exhibited on triangular cells. In particular, for general
constructions there may be loss of consistency for the $\mathbb{V}_1$
spaces used here (the consistency problem for $\mathbb{V}_2$ is
avoided through rehabilitation). However, \citet{Holst12} demonstrated
that consistent approximation can still be obtained if the
computational domain can be obtained \emph{via} a transformation from
a mesh of affine elements onto the higher-order nodal interpolation of
a $C^\infty$ manifold (such as the sphere). This aspect is discussed
further in the context of geophysical fluid dynamics in
\cite{natale2016compat}.

\subsubsection{Mixed finite element formulation}

The mixed finite element discretisation of equations
\eqref{eq:momentum}-\eqref{eq:mass} is formed by restricting $\MM{u}\in
\mathbb{V}_1$, $D\in \mathbb{V}_2$, multiplying the equations by
appropriate test functions, $\MM{w}\in\mathbb{V}_1$ and
$\phi\in\mathbb{V}_2$, and integrating over the domain, giving
\begin{align}
\int_\Omega\MM{w}\cdot\MM{u}_t \diff V + \int_\Omega\MM{w}\cdot \MM{Q}^\perp \diff V - \int_\Omega\nabla\cdot\MM{w} \left(g\left(D+b\right)+\frac{1}{2}|\MM{u}|^2\right) \diff V &= 0, \quad \forall \MM{w}\in \mathbb{V}_1,
\label{eq:weak-momentum}\\
\int_\Omega \phi\left(D_t + \tilde{\nabla}\cdot\MM{F}\right) \diff V &= 0, 
\quad \forall \phi \in \mathbb{V}_2. \label{eq:weak-mass}
\end{align}
where we have introduced the mass flux $\MM{F}\in \mathbb{V}_1$ (an
approximation to $\MM{u}D$), and the vorticity flux $\MM{Q}$ (an
approximation to $\MM{u}Dq$). $\MM{F}$ and $\MM{Q}$ are defined by
appropriate choice of advection schemes for $D$ and the diagnostic
potential vorticity $q$, which we shall define later.  

We have integrated the gradient term in equation \eqref{eq:weak-momentum}
by parts to avoid taking the gradient of the layer depth which is
undefined since $D\in\mathbb{V}_2$ is discontinuous. A similar 
problem is posed for the definition of the potential vorticity $q$
since $\zeta=\nabla^\perp\cdot\MM{u}$, which appears in the definition of
$q$, is similarly undefined. In order to fix this we integrate the
$\nabla^\perp$ term by parts (with no surface term since we are
solving our equations on the sphere), defining our discrete
potential vorticity $q\in\mathbb{V}_0$ as
\begin{equation}
\label{eq:discretePV}
\int_\Omega \gamma qD \diff V = \int_\Omega -\nabla^\perp\gamma\cdot\MM{u}\diff V + \int_\Omega\gamma f\diff V \quad \forall\gamma\in\mathbb{V}_0.
\end{equation}
Provided $D$ remains positive (this can be enforced using a slope
limiter, although in the test cases used here the mean depth is
sufficiently high that a slope limiter is not required), then this
equation can be solved for $q\in \mathbb{V}_0$ from known $D$ and
$\MM{u}$. This provides a one-to-one mapping between $q$ and the weak
curl of $\MM{u}$, which leads to a one-to-one mapping between $q$ and
the divergence-free part of $\MM{u}$ \emph{via} the discrete Helmholtz
decomposition (without needing to solve a global Poisson
problem). Hence, control of $q$ in the $L^2$ norm from a chosen stable
advection scheme provides strong control over the divergence-free
component of $\MM{u}$ without compromising the divergent part. This is
one reason for the success of this type of potential vorticity
conserving discretisation for the rotating shallow water equations.

We differentiate this equation in time, and substitute for $\MM{u}_t$
using equation \eqref{eq:weak-momentum} with
$\MM{w}=-\nabla^\perp\gamma$. Since $\nabla\cdot\nabla^\perp\equiv0$
and we assume $f_t=0$, this gives
\begin{equation}
\label{eq:PV-conservation-inV0}
\int_\Omega\gamma(qD)_t\diff V - \int_\Omega\nabla\gamma\cdot \MM{Q}\diff V=0,
\quad \forall \gamma \in \mathbb{V}_0,
\end{equation}
which is the Galerkin projection of the PV conservation law into
$\mathbb{V}_0$. If we select $\gamma=1$, then we obtain conservation
of total PV,
\begin{equation}\label{eq: int qD}
\dd{}{t}\int_\Omega qD \diff V = 0.
\end{equation}
On the other hand, if we are on the sphere, then this quantity is zero
as can be computed directly from \eqref{eq:discretePV}; this
topological result stems from the fact that in this case
$\MM{w}=-\nabla{^\perp}\gamma=0$.

If we choose our advection scheme so that $\MM{Q}=q\MM{F}$
{and} if $q$ is a constant, then we may integrate
{\eqref{eq:PV-conservation-inV0}} by parts without introducing error
(since $q\in H^1$ and $\MM{F}\in H(\mbox{div})$), and we obtain
\begin{equation}
\int_\Omega\gamma(q D)_t\diff V = -\int_\Omega \gamma {q}\nabla\cdot(\MM{F})\diff V,
\end{equation}
and hence
\begin{equation}
\int_\Omega\gamma D q_t\diff V = -\int_\Omega \gamma q\left(D_t + \nabla\cdot\MM{F}
\right)\diff V.
\end{equation}
For flat elements, or if we had chosen Option 1, the right hand side
is zero since $D$ and $\nabla\cdot\MM{F}$ are in the same finite
element space and therefore equation \eqref{eq:weak-mass} holds pointwise,
\emph{i.e.},
\begin{equation}
\label{eq:D pwise}
D_t + \nabla\cdot\MM{F}=0.
\end{equation}
Hence, we conclude that if $q$ is constant, then $q_t=0$, and so $q$
remains constant. As well as being a statement of first-order consistency
for the advection scheme for $q$, this is also an important property
of equation \eqref{eq:pv advection}.

The formulation requires some adaptation for the finite element spaces
used in this paper to obtain this result, since $\nabla\cdot\MM{F}$ is
not guaranteed to be in the same space as $D$, so equation \eqref{eq:D
  pwise} does not hold pointwise. To recover it, some further
``rehabiliation'' is required; we amend equation \eqref{eq:discretePV} by
defining $\tilde{D}\in\correction{\mathbb{V}_2}$ such that
\begin{equation}
\frac{\tilde{D}_t}{\tau}+\nabla\cdot\MM{F}=0,
\end{equation}
where
\begin{equation}
\tau|_e = \det J_e\circ g_e.
\end{equation}
Comparing the weak form of this equation with \eqref{eq:weak-mass} we see
that
\begin{equation}
\int_\Omega\phi\frac{\tilde{D}_t}{\tau}\diff V = \int_\Omega\phi D_t\diff V \quad \forall\phi\in\mathbb{V}_2,
\end{equation}

\noindent which is consistent with the definition

\begin{equation}
\int_\Omega\phi\frac{\tilde{D}}{\tau}\diff V = \int_\Omega \phi D\diff V \quad \forall\phi\in\mathbb{V}_2,
\end{equation}
hence we can solve for $\tilde{D}$ (since $\tau$ is a positive
quantity and just alters the metric). $\tilde{D}$ is discontinuous and
hence this equation can be solved separately in each element.

Using this in equation
\eqref{eq:discretePV} gives

\begin{equation}
\label{eq:PVadvection}
\int_\Omega\gamma q\frac{\tilde{D}}{\tau} \diff V = -\int_\Omega\nabla^\perp\gamma\cdot\MM{u} \diff V + \int_{\Omega}\gamma f \diff V \quad \forall\gamma\in\mathbb{V}_0.
\end{equation}

\noindent Differentiating, rearranging and assuming $q$ is constant as
before we obtain
\begin{equation}
\int_\Omega\gamma\frac{\tilde{D}}{\tau}q_t\diff V = -\int_\Omega\gamma q\left(\frac{\tilde{D}{_t}}{\tau}+\nabla\cdot\MM{F}\right)\diff V = 0,
\end{equation}
as required.

In fact, the finite element formulation allows us to go beyond the
first order consistency result. For example, if we choose
$\MM{Q}=\MM{F}q$, we have enough continuity to integrate by parts, 
and we obtain
\begin{equation}
\int_\Omega \gamma\left((qD)_t + \nabla\cdot(\MM{F}q)\right)\diff V = 0,
\quad \forall \gamma \in \mathbb{V}_0.
\end{equation}
The left-hand side vanishes if $q$ is an exact solution of the equation
\begin{equation}
\label{eq:exact}
(qD)_t + \nabla\cdot(\MM{F}q)=0,
\end{equation}
where $D$ and $\MM{F}$ are the discrete mass and mass flux, indicating
that the discretisation is consistent at the order of approximation of
the finite element space $\mathbb{V}_0$. Unfortunately, this is not a
good choice for $\MM{Q}$ since the implied discrete advection operator
is not stable. An alternative is to choose
\begin{equation}
\MM{Q} = \MM{F}q - \alpha\MM{F}\frac{\MM{F}}{|\MM{F}|}\cdot\nabla q,
\end{equation}
where $\alpha$ is a stabilisation parameter. Substitution, rearrangement
and integration by parts on the $\MM{F}q$ term leads to
\begin{equation}
\int_\Omega \bigg(\gamma + \alpha\frac{\MM{F}}{|\MM{F}|}\cdot\nabla{\gamma}\bigg)\big((qD)_t + \nabla\cdot(\MM{F}q)\big)\diff V = 0, 
\quad \forall \gamma \in \mathbb{V}_0,
\end{equation}
which is a streamline-upwind Petrov-Galerkin (SUPG) spatial
discretisation of the PV conservation equation. This discretisation is
stable, and also consistent, \emph{i.e.} the equation vanishes when
the exact solution to \eqref{eq:exact} is substituted.

In this paper we will use a Taylor-Galerkin discretisation for the
potential vorticity equation; this discretisation achieves the same
aims as the SUPG discretisation described above, but arises more
naturally in the discrete time setting and hence we shall postpone our
discussion of it until we have described the time-discrete formulation
of the full shallow water system in the next section.

{We remark that all of the results above are assuming that
  integrations are evaluated exactly. Due to the factors of $1/\det J$
  arising from the Piola transformation, not all integrands are
  polynomial and hence exact quadrature is difficult. However, as
  noted in \cite{cotter2014finite}, the structure of terms of the
  following forms,
\begin{equation}
  \int_\Omega \nabla g \cdot \MM{v} \diff x, \quad
  \int_\Omega g \MM{w}\cdot \MM{v}^\perp \diff x, \quad
  \quad \int_\Omega g \nabla\cdot\MM{w} \diff x,
\end{equation}
where $g$ is an arbitrary scalar function (such as a product of scalar
functions) and $\MM{w}$ and $\MM{v}$ are Piola-mapped functions,
means that the factors of $\det J$ cancel after transforming to
the reference element, and the integral after change of variables
has a polynomial integrand after all. Later on this will also apply
to integrals of the form
\begin{equation}
  \int_f \MM{v}\cdot\MM{n} g \diff S,
\end{equation}
where $f$ is an element facet. This means that all of the conservation
properties in this paper also hold for inexact quadrature provided
that it is sufficiently high order to integrate these polynomial
integrands exactly, after appropriately redefining inner products
using a quadrature rule instead of an exact integral.}

\subsection{{Time discretisation}}

We shall build a semi-implicit time discrete formulation. First we 
write
\begin{equation}
\MM{u}^*=\theta\MM{u}^{n+1}+(1-\theta)\MM{u}^n, D^*=\theta D^{n+1} +(1-\theta)D^n
\end{equation}
and
\begin{equation}
\Delta\MM{u}=\MM{u}^{n+1}-\MM{u}^n, \Delta D=D^{n+1}-D^n.
\end{equation}
We can now write equations \eqref{eq:weak-momentum} and \eqref{eq:weak-mass} as
\begin{align}
\label{eq:du/dt discrete time}
\int_\Omega\MM{w}\cdot\Delta\MM{u}\diff V
+\Delta t\int_\Omega\MM{w}\cdot\MM{Q}^\perp\diff V -\Delta t\int_\Omega\nabla\cdot\MM{w}\left(g({D}^*+b)+K(\MM{u}^*)\right)\diff V &= 0,
\quad \forall \MM{w}\in \mathbb{V}_1, \\
\int_\Omega\phi\Delta D\diff V+\Delta t\int_\Omega\phi\nabla\cdot\MM{F}\diff V=&0,
\quad \forall \phi \in \mathbb{V}_2,
\end{align}
where 
\begin{equation}
K(\MM{u}) = \frac{1}{2}|\MM{u}|^2,
\end{equation}
and where the time-averaged mass and vorticity fluxes $\MM{F}$ and
$\MM{Q}$ are yet to be defined. The idea is that we choose $\MM{F}$ to
be a time-independent function such that $D^{n+1}$ is obtained from
$D^n$ \emph{via} a high-order stable time discretisation over one
timestep for the equation
\begin{equation}
D_t + \nabla\cdot(\MM{u}^*D) = 0,
\end{equation}
\emph{i.e.} with the advecting velocity frozen to the value of
$\MM{u}^*$. Similarly, $\MM{Q}$ is chosen so that $q^{n+1}$ is related
to $q^n$ \emph{via} a high-order stable time discretisation over one
timestep for the equation
\begin{equation}
(qD)_t + \nabla\cdot(\MM{F}q)=0,
\end{equation}
\emph{i.e.} with the advecting mass flux frozen to the time-averaged
flux $\MM{F}$. This means that for $\theta=1/2$ we obtain a scheme
that is overall second-order in time, but that uses higher-order
advection schemes for $D$ and $q$. The rationale is that for
near-linear waves, we would like the propagation to be as conservative
as possible, so the semi-linear formulation should be based around a
time-centred scheme. However, we would also like to obtain good
quality solutions over long integrations when close to geostrophic
balance, in which case the important quantity is $q$, and it is
important that we transport $D$ consistently with $q$ to stay close to
the balanced state. In that regime, it is thought to be important to
use a high odd-order time integration scheme, since for odd-order
schemes the error is dominated by diffusion rather than dispersion
(the latter leads to oscillations near to near-discontinuous data).
It
is also thought that the use of a time-averaged velocity to transport
$q$ and $D$ helps to preserve geostrophic balance.
This
was also the rationale behind choosing the 3rd order Forward-in-Time
advection schemes used in \citet{thuburn2015primal}.

An implicit formulation such as the one above requires Newton or
Picard iterations to iterate to convergence. In practice, we perform a
small fixed number {(4)} of Picard iterations per timestep,
since our aim is to obtain a stable timestepping method with accurate
$q$ and $D$ transport, rather than iterating to convergence to obtain
an exact implementation of the fully implicit scheme described
above. In the Picard iteration we replace the Jacobian obtained from
linearisation around the current state with the Jacobian linearised
around the system at a state of rest. This means that it is possible
to reduce the system as described in the next section, and that the
same solver context can be reused during each Picard iteration and
timestep.

Hence, each Picard cycle consists of the following steps.
\begin{enumerate}
\item Initialise $\Delta \MM{u}=0$, $\Delta D=0$.
\item Use the current values of $\Delta \MM{u}$ and $\Delta D$ to
  compute corresponding values for $\MM{u}^*$, $D^*$.
\item Use the chosen mass advection scheme with velocity $\MM{u}^*$
  to update from $D^n$ to $D^{n+1}$, and find $\MM{F}$ such that
\begin{equation}
\int_\Omega \phi\left({D}^{n+1}-{D}^n + \Delta t\nabla\cdot\MM{F}\right)\diff V {= 0},
\quad \forall \phi\in \mathbb{V}_2.
\end{equation}
We will explain how to construct $\MM{F}$ in section \ref{sec:mass adv}.
The mass residual $R_D:\mathbb{V}_2\to \mathbb{R}$ is defined as
\begin{equation}
R_D[\phi] = \int_\Omega \phi(\Delta D + \Delta t\nabla\cdot\MM{F}{)}\diff V.
\end{equation}
\item Diagnose the PV $q^n$ at time $t^n$. 
\item Use the chosen PV advection scheme with mass flux $\MM{F}$ 
to update from $q^n$ to $q^{n+1}$ and compute the corresponding PV
flux $\MM{Q}$ such that
\begin{equation}
\int_\Omega \gamma \left(q^{n+1}\frac{\tilde{D}^{n+1}}{\tau}
- q^n\frac{\tilde{D}^n}{\tau}\right)\diff V
- \Delta t\int_\Omega \nabla\gamma \cdot \MM{Q}\diff V = 0,
\quad \forall \gamma\in \mathbb{V}_0.
\end{equation}
\item The velocity residual $R_{\MM{u}}:\mathbb{V}_1\to\mathbb{R}$ is
  defined as
\begin{equation}
  R_{\MM{u}}[\MM{w}] = \int_\Omega\MM{w}\cdot\Delta\MM{u}\diff V
  +\Delta t\int_\Omega\MM{w}\cdot\MM{Q}{^\perp}\diff V +\Delta t\int_\Omega\nabla\cdot\MM{w}\left(g({D}^*+b)+K(\MM{u}^*)\right)\diff V.
\end{equation}
\item The increments 
\begin{equation}
\Delta \MM{u} \mapsto \Delta \MM{u} + \delta \MM{u}, 
\quad \Delta D \mapsto \Delta D + \delta D
\end{equation}
are then obtained by solving the coupled system
\begin{align}
\label{eq:delta u}
\int_\Omega \MM{w}\cdot\delta\MM{u}\diff V + \theta\Delta t\int_\Omega
f\MM{w}\cdot\delta\MM{u}^\perp \diff V - \theta\Delta t\int_\Omega \nabla\cdot\MM{w}
g\delta D \diff V & = -R_{\MM{u}}[\MM{w}], \quad \forall \MM{w}\in \mathbb{V}_1, \\
\label{eq:delta D}
\int_\Omega \phi\left(\delta D + \theta\Delta t H_0\nabla\cdot\delta\MM{u}\right)
\diff V & = -R_D[\phi], \quad \forall \phi\in\mathbb{V}_2,
\end{align}
where $H_0$ is the mean layer depth at rest.
\item If {we have not completed 4 iterations}, apply these updates and
return to 2.
\end{enumerate}

In the following sections we describe the construction of
$\MM{F}$, $\MM{Q}$ and the solution of the coupled system in detail.

\subsubsection{Solving the coupled linear system}

In our formulation, we use hybridisation to solve equations
(\eqref{eq:delta u}-\eqref{eq:delta D}). Hybridisation is a technique for
efficiently solving mixed finite element problems that has been used
since the 1960s; in the 1980s it was also discovered that the
hybridised formulation could also be used to obtain more accurate
approximations of the solution (see \cite{boffi2013mixed} for a
general survey).

Obtaining a hybridised formulation requires two steps. First, we
introduce a finite element space $\tilde{\mathbb{V}}_1$ by relaxing
the normal continuity constraints within $\mathbb{V}_1$. In other
words, functions in $\tilde{\mathbb{V}}_1$ have the same local
polynomial representation as functions in $\mathbb{V}_1$, but there
are no requirements of continuity between edges. In particular, we
note that $\mathbb{V}_1\subset \tilde{\mathbb{V}}_1$. Second, we
introduce a trace space $\Tr(\mathbb{V}_1)$, defined on the element
facet set $\Gamma$, such that functions $\lambda \in
\Tr(\mathbb{V}_1)$ are scalar functions which when restricted to a
single element facet $f$, are from the same polynomial space as
$\MM{u}\cdot\MM{n}$ restricted to that facet. Having relaxed the
continuity requirements for $\delta\MM{u}\in\tilde{\mathbb{V}}_1$,
we enforce them again by adding another equation,
\begin{equation}
\label{eq:constraint}
\int_\Gamma \mu\llbracket \MM{u} \rrbracket dS = 0, \quad \forall \mu \in
\Tr(\mathbb{V}_1),
\end{equation}
where we use the usual ``jump'' notation
\begin{equation}
\llbracket \MM{u} \rrbracket = \MM{u}^+\cdot\MM{n}^+ + \MM{u}^-\cdot\MM{n}^-,
\end{equation}
having arbitrarily labelled each side of each facet with $+$ and $-$,
so that $\MM{n}^+$ points from the $+$ side to the $-$ side and vice
versa. To avoid an over-determined system, we introduce Lagrange multipliers
$\lambda\in \Tr(\mathbb{V}_1)$ and rewrite equations (\eqref{eq:delta u}-\eqref{eq:delta D}) as
\begin{align}
\label{eq:delta u hybridised}
\int_\Omega \MM{w}\cdot\delta\MM{u}\diff V + \theta\Delta t\int_\Omega
f\MM{w}\cdot\MM{u}^\perp \diff V - \theta\Delta t\int_\Omega \nabla\cdot\MM{w}
g\delta D \diff V + 
\int_\Gamma \lambda \llbracket\MM{w}\rrbracket \diff S
& = -R_{\MM{u}}[\MM{w}], \quad \forall \MM{w}\in \tilde{\mathbb{V}}_1, \\
\int_\Omega \phi\left(\delta D + \theta\Delta t H_0\nabla\cdot\delta\MM{u}\right)
\diff V & = -R_D[\phi], \quad \forall \phi\in\mathbb{V}_2,
\label{eq:delta D hybridised}
\end{align}
together with equation \eqref{eq:constraint}. Note that the residual
$R_{\MM{u}}$ must now be evaluated with $\MM{w}\in
\tilde{\mathbb{V}}_1$.  All of the inter-element coupling in equations
(\eqref{eq:delta u hybridised}-\eqref{eq:delta D hybridised}) takes place
in the $\lambda$ term. This means that if $\lambda$ is known, then it
is possible to obtain $\MM{u}$ and $D$ independently in each
element. To enable elimination of $\MM{u}$ and $D$, we define a
lifting operator $L:\Tr(\mathbb{V}_1)\to \tilde{\mathbb{V}}_1)$, which
gives the solution $\MM{u}$ for a given $\lambda$ in the case when
$R_D[\phi]$ and $R_{\MM{u}}[\MM{w}]$ are replaced by zero. We also
define $\MM{u}_0$ which is the solution obtained when $\lambda$ is
zero, but $R_D[\phi]$ and $R_{\MM{u}}[\MM{w}]$ are present. Then,
the general solution of this equation given particular $R_D$ and 
$R_{\MM{u}}$ is
\begin{equation}
\MM{u} = L\lambda + \MM{u}_0.
\end{equation}
Substituting into equation \eqref{eq:constraint} gives
\begin{equation}
\int_\Gamma \mu\llbracket L\lambda \rrbracket \diff S = 
-\int_\Gamma \mu\llbracket \MM{u}_0 \rrbracket \diff S,
\quad \forall \mu \in \Tr(\mathbb{V}_1).
\label{eq:hybridised}
\end{equation}
This reduced equation can be solved for $\lambda$ before
reconstructing $\MM{u}$ and $D$ by solving equations (\eqref{eq:delta
  u hybridised}-\eqref{eq:delta D hybridised}) independently in each
element. Since the value of $L\lambda$ in each element only depends on
the values of $\lambda$ on the facets of that element, equation
\eqref{eq:hybridised} only couples together values of $\lambda$ on
facets that share an element. This means that the matrix-vector form
of this equation is very sparse. In fact, the matrix can be assembled
by visiting each element separately, performing inversion on element
block systems. Further, this equation has the same spectral properties
as the Helmholtz operator, and hence can be solved with Krylov methods
and standard preconditioners such as SOR, algebraic multigrid; a
geometric multigrid for general higher-order hybridised mixed finite
element elliptic problems was provided in
\cite{gopalakrishnan2009convergent}. Note that the Coriolis term can
be included in the linear system in this approach without altering the
sparsity of the reduced system. One important aspect is that if the
solver for this system is not iterated to convergence, the resulting
velocity field will not be exactly {div-conforming}. We
address this in our implementation by simply projecting the velocity
back into {$\mathbb{V}_1$} after reconstruction, which
appears not to cause any problems with stability. A more sophisticated
approach would use the hybridised solver as a preconditioner for the
$(\Delta u, \Delta {D})$ system on
{$(\mathbb{V}_1, \mathbb{V}_2)$}; we will investigate this
in further work.

\subsubsection{Advection scheme for layer depth $D$}
\label{sec:mass adv}
We now need to solve the mass continuity equation \eqref{eq:weak-mass} for
the update to $D$. As $D\in \mathbb{V}_2$ and is therefore
discontinuous, we can use standard upwind discontinuous Galerkin methods
to obtain an approximation to $D_t$,
\begin{equation}
\label{eq:DG Euler}
\int_e\phi D_t\diff V=\Delta t\int_e\nabla\phi\cdot\MM{u}^*{D}\diff V - \Delta t\int_{\partial e}\phi D^u\MM{u}^*\cdot\MM{n} \diff S, \quad \forall \phi \in \mathbb{V}_2(e),
\end{equation}
for each element $e$, where $\mathbb{V}_1(e)$ is the space
$\mathbb{V}_1$ restricted to the element $e$, and $D^u$ is the value of
$D$ on the upwind side of the element boundary $\partial e$.  We then
use the standard 3-stage Strong Stability Preserving Runge-Kutta
scheme \citep{gottlieb2001strong},
\begin{align}
\label{eq:ssprk3-1}
D^1 & = D^n + \Delta t D_t^n, \\
D^2 & = \frac{3}{4}D^n + \frac{1}{4}\left(
D^1 + \Delta t D_t^1
\right), \\
D^{n+1} & = \frac{1}{3}D^n + \frac{2}{3}\left(
D^2 + \Delta t D_t^2
\right).
\label{eq:ssprk3-3}
\end{align}
Later we shall discuss a consistency property of the potential
vorticity conserving discretisation of the velocity equation; this property requires
that we find a time-integrated mass flux $\bar{\MM{F}}\in\mathbb{V}_1$
that satisfies
\begin{equation}
\label{eq:deltaD}
\Delta D = -\Delta t\nabla\cdot\bar{\MM{F}}.
\end{equation}
This is done by finding $\MM{F}$ given $D_t$ such that
\begin{equation}
\label{eq:discrete cty}
D_t + \nabla\cdot\MM{F} = 0,
\end{equation}
and then substituting into equations (\eqref{eq:ssprk3-1}-\eqref{eq:ssprk3-3})
to construct $\bar{\MM{F}}$. To find a flux $\MM{F}$ for each 
Runge-Kutta stage, we solve
\begin{align}
\int_{f}\kappa\MM{F}\cdot\MM{n} \diff s &= \int_{f}\kappa D^u\MM{u}^*\cdot\MM{n}\diff s & \forall \kappa \in \Tr{_f}(\mathbb{V}_1), {\quad \forall f \in \partial e}\label{eq:Fdotn}\\
\int_e\MM{w}\cdot\MM{F}\diff x &= \int_e\MM{w}\cdot\MM{u}^*D\diff x & \forall \phi \in \mathring{\mathbb{V}}_1^*(e), \label{eq:Fdx}
\end{align}
where $\mathring{\mathbb{V}}^*_1(e)$ is the right size to close the
system and $\nabla\phi \in \mathring{{\mathbb{V}}}_1^*(e)$ for all
$\phi \in \mathbb{V}_2(e)$. The left-hand sides of equations
(\eqref{eq:Fdotn}-\eqref{eq:Fdx}) are the same as the left-hand sides
of the definition of the commuting operator $\pi:\Hdiv\to
\mathbb{V}_1$ that features in stability proofs for mixed finite
element methods (see \cite{boffi2013mixed}, for example). This means
that the above construction is well-posed. To check that equation
\eqref{eq:discrete cty} is satisfied, we integrate \eqref{eq:deltaD} by
parts. Then, substituting the above relations
\eqref{eq:Fdotn}-\eqref{eq:Fdx} (with $\MM{w}=\nabla\phi)$, we see
that

\begin{align*}
\int_e\phi D_t\diff V &= -\Delta t\int_e\phi\nabla\cdot\MM{F}\diff V \\
                           &=  \Delta t\int_e\nabla\phi\cdot\MM{F}\diff V - \Delta t\int_{\delta e}\phi\MM{F}\cdot\MM{n}\diff S \\
                           &=  \Delta t\int_e\nabla\phi\cdot\MM{u}^*D\diff V - \Delta t\int_{\delta e}\phi D^u\MM{u}^*\cdot\MM{n}\diff S,
\end{align*}

\noindent as required. 

We note, although we did not use it in this paper, that the use of a
slope limiter can also be incorporated into the mass flux computation.
If a slope limiter is used after a Runge-Kutta stage, then $D$ is
replaced by $D'$, with $\int_e D-D'\diff x =0$. If we seek $\MM{F}'$ 
such that
\begin{equation}
D'-D + \Delta t \nabla \cdot \MM{F}' = 0,
\end{equation}
then integration over a single element immediately tells us that 
this can be satisfied by
\begin{equation}
\int_f \kappa \MM{F}'\cdot\MM{n} \diff S = 0, \quad 
\forall \kappa \in \Tr{_f}(\mathbb{V}_1), {\quad \forall f \in \partial e},
\end{equation}
\emph{i.e.} $\MM{F}'\cdot\MM{n}=0$ on the boundary $\partial e$. We
then solve a local mixed problem for $(\MM{F}',p) \in
(\mathbb{V}_1(e), \mathbb{V}_2(e))$, given by
\begin{align}
\int_e \phi \nabla\cdot\MM{F}' \diff x &= \int_e \phi(D'-D)\diff x, {\quad \forall \phi \in \mathbb{V}_2(e)},\\
\int_e \MM{w}\cdot\MM{F}' + \nabla\cdot\MM{w}p \diff x & = 0, {\quad \forall \MM{w} \in \mathbb{V}_1(e)},
\end{align}
subject to the above zero Dirichlet boundary conditions for $\MM{F}'$,
{where $p$ is introduced to determine a unique $\MM{F}'$}.


\subsubsection{Advection scheme for velocity $\MM{u}$}

The advection scheme described in this section follows the following
design strategy: find an advection scheme for $q$ that is compatible
with equation \eqref{eq:du/dt discrete time}, and select the
corresponding $\MM{Q}$ for insertion into that equation to compute
$R_{\MM{u}}[\MM{w}]$. Selecting $\MM{w}=-\nabla^\perp\gamma$ in equation
\eqref{eq:du/dt discrete time}, evaluating equation \eqref{eq:discretePV}
at time levels $n$ and $n+1$ and substituting gives 
\begin{equation}
\label{eq:discrete q}
\int_\Omega \gamma q^{n+1}D^{n+1} \diff x 
-\int_\Omega \gamma q^{n}D^{n} \diff x 
- \Delta t \int_\Omega\nabla \gamma \cdot\MM{Q}\diff x = 0, 
\quad \forall \gamma \in \mathbb{V}^0,
\end{equation}
after noting that $\MM{w}$ is divergence-free.
\citet{mcrae2014energy} used $\MM{Q}=\bar{\MM{F}}q^{n+1/2}$ to obtain
an implicit midpoint rule time discretisation for an energy-enstrophy
conserving formulation. In this paper, we aim to use higher-order
stabilised advection schemes in order to obtain accurate
representation of potential vorticity transport. We note that
Streamline Upwind Petrov Galerkin methods \citep{brooks1982streamline}
and Taylor-Galerkin methods \citep{donea1984taylor} all result in
time-discrete formulations equivalent to equation \eqref{eq:discrete
  q}.

It is desirable to use a higher order timestepping scheme for $q$, to
obtain accurate transport of potential vorticity using the balanced
flow. In particular, odd-ordered schemes are attractive since the
leading order error is diffusive rather than dispersive. Hence, in
this paper we make use of the two-step third-order unconditionally-stable
Taylor-Galerkin scheme of \citet{safjan1993high}.

Taylor-Galerkin schemes are built by transforming time derivatives into
space derivatives using the advection equation. 

The general form of a multistage Taylor-Galerkin method is
\begin{equation}
  \label{eq:TG}
{Z}_i - \eta(\Delta t)^2(Z_i)_{tt} = Z^n + \Delta t\sum_{j=1}^i\mu_{ij}(Z_j)_t + (\Delta t)^2\sum_{j=1}^i\nu_{ij}(Z_j)_{tt}, \quad i=1,\ldots,k,
\end{equation}
where $\eta$ is a stabilisation parameter, the subscript $i$
is the stage index, and the $\{\mu\}_{ij}$ and $\{\nu\}_{ij}$
are coefficients defined in \citet{safjan1993high}. After computing
these stage variables, the value of $Z_k$ is copied into $Z^{n+1}$.

In our discretisation we use a formulation where $Z=qD$. Then we
have
\begin{align}
  \int_\Omega \gamma (qD)_t \diff x & = \int_\Omega \nabla\gamma\cdot
  \left(\bar{\MM{F}}q\right)\diff x, \quad \forall\gamma\in \mathbb{V}_0, \\
  \int_\Omega \gamma (qD)_{tt} \diff x & = \int_\Omega \nabla\gamma\cdot
  \left(\bar{\MM{F}}q_t\right)\diff x
  = -\int_\Omega \nabla\gamma\cdot
  \left(\bar{\MM{F}}\frac{\bar{\MM{F}}}{\bar{D}}\cdot {\nabla} q\right)\diff x, 
\quad \forall\gamma\in \mathbb{V}_0,
\end{align}
recalling that $\bar{F}=\bar{D}\bar{\MM{u}}$ is considered to be
time-independent over the advection step as part of the semi-implicit
discretisation. Combination with Equation \eqref{eq:TG} leads to
\begin{align}
\nonumber
\int_\Omega \gamma q_iD^{n+1} + \eta\Delta t^2\frac{\bar{\MM{F}}}{\bar{D}}\cdot\nabla\gamma\bar{\MM{F}}\cdot\nabla q_i\diff x & =
\int_\Omega \gamma q^nD^n + \Delta t\sum_{j=1}^i\mu_{ij}\int_\Omega\nabla\gamma\cdot\bar{\MM{F}}q_j \diff x - \\
& \qquad (\Delta t)^2\sum_{j=1}^i\nu_{ij}
\int_\Omega\frac{\bar{\MM{F}}}{\bar{D}}\cdot\nabla\gamma\bar{\MM{F}}\cdot\nabla q_j\diff x,
\quad \forall \gamma \in \mathbb{V}_0.
\label{eq:multistage}
\end{align}
Note that this equation involves solving a Helmholtz-type equation
with derivatives in the streamline direction for each stage variable.
This equation is symmetric positive definite and well-conditioned for
$\mathcal{O}(1)$ Courant numbers; the conjugate gradient method
converges quickly with simple SOR preconditioning.

In this paper we took the following coefficients
\begin{align}
\nonumber
\eta = 0.48,\,
c_1 = \frac{1}{2}\left(1 + (-\frac{1}{3}+8\eta)^\frac{1}{2}\right),\,
\mu_{11} = c_1,\,
\mu_{12} = 0, \,
\mu_{21} = \frac{1}{2}\left(3-\frac{1}{c_1}\right),\,\\
\mu_{22} = \frac{1}{2}\left(\frac{1}{c_1}-1\right),\,
\nu_{11} = \frac{1}{2}c_1^2-\eta,\,
\nu_{12} = 0,\,
\nu_{21} = \frac{1}{4}\left(3c_1-1\right)-\eta,\,
\nu_{22} = \frac{1}{4}\left(1-c_1\right).
\end{align}
\citet{safjan1993high} showed this scheme to be third order and
unconditionally stable for $\eta>0.473$.

Having solved for $q^{n+1}$, we take $i=2$, and notice that equation
\eqref{eq:multistage} takes the form of equation \eqref{eq:discrete
  q}, and hence $\MM{Q}$ can be extracted for insertion into equation
\eqref{eq:du/dt discrete time}. For the case of curved elements, $D$
must be replaced by $\tilde{D}/\tau$ throughout.

\section{Numerical Results}
\label{sec: results}

In this section we show numerical results from three standard test
cases on icosahedral grids, using the spaces $(P_3,BDM_2,DG_1)$, and a
piecewise cubic approximation to the surface of the sphere. The code
was implemented using the {Firedrake finite element
  framework, which permits the symbolic implementation of the mixed
  finite element techniques discussed in this paper. Using the Unified
  Form Language (see \citet{UFL}), which provides a high-level
  mathematical description of the finite element problem,
  code is automatically generated which forms the resulting matrix equations
  by assembling the local contributions from each cell and/or facet of the
  mesh. These equations are provided directly to PETSc \citep{petsc-user-ref, petsc-efficient},
  which provides direct access to runtime configurable iterative solvers
  and preconditioners. The hybridisation of the implicit system for the linear
  updates is implemented using the Slate framework \citep{gibson2018domain},
  which performs the local elimination and recovery operations. The reduced equation
  for the trace variables is numerically inverted using the conjugate gradient
  method and PETSc's smooth aggregation multigrid preconditioner (GAMG).}

\begin{figure}
  \centering
  \includegraphics[width=0.25\linewidth]{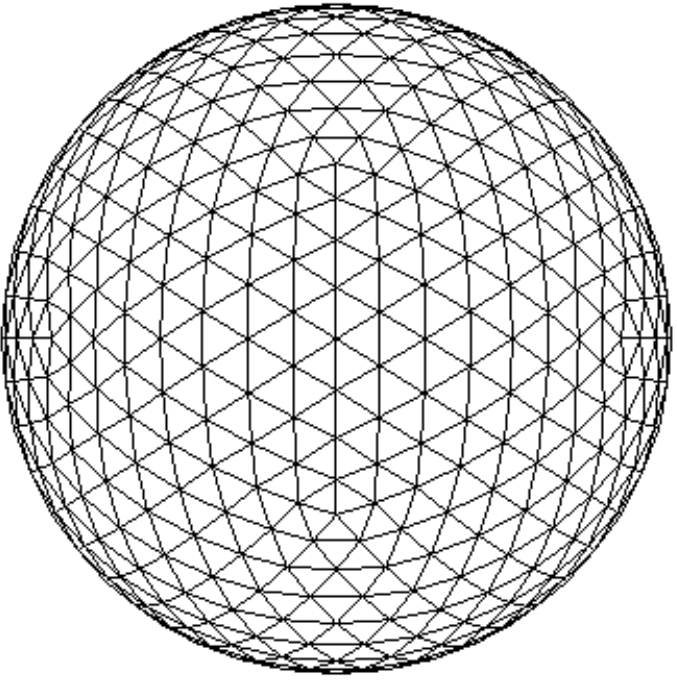}
  \caption{Icosahedral sphere grid, viewed looking down on the North pole, corresponding to the lowest resolution described in table \ref{table:grid params}.}
  \label{fig:sphere_grid}
\end{figure}

The first two test cases are numbers $2$ (solid body rotation) and $5$
(flow over a mountain) from \citet{williamson1992standard} and the
third is the barotropically unstable jet from
\citet{galewsky2004initial}. Table \ref{table:grid params} contains
information on the properties of the 4 grids that we use for the
Williamson convergence tests, along with the timesteps used. The
timestep was chosen to give a constant Courant number across the
different resolutions; $0.2$ in the solid body rotation test and
$0.06$ in the flow over a mountain test, {although we note
  that in order to see second order convergence, we had to further
  reduce the timestep for the highest resolution mountain simulation
  because at this resolution the time discretisation errors start to
  dominate.} Figure \ref{fig:sphere_grid} shows the lowest resolution
icosahedral grid that we use.

\begin{table}
  \centering
  \begin{tabular}{@{\extracolsep{6pt}}rrrrrrrrS}
    \hline
    \multicolumn{4}{c}{Grid properties} & \multicolumn{3}{c}{DOFs} & {Test 2} & {Test 5} {\bigstrut} \\\cline{1-4}\cline{5-7}\cline{8-8}\cline{9-9}
  {cells} & {nodes} & {$\Delta x_{\max}$ (km)} & {$\Delta x_{\min}$ (km)} & {$\MM{u}$} & {$D$} & {$q$} & {$\Delta t$ (s)}  & {$\Delta t$ (s)}  {\bigstrut} \\
  \hline
  1280 & 642 & 1054 & 720 & 9600 & 3840 & 5762 & 3000 &  900 {\bigstrut[t]} \\
  5120 & 2562 & 527 & 348 & 38400 & 15360 & 23042 & 1500 & 450 \\
  20480 & 10242 & 263 & 171 & 153600 & 61440 & 92162 & 750 & 225 \\
  81920 & 40962 & 132 & 85 & 614400 & 245760 & 368642 & 375 & 84.375 \\
  \end{tabular}
  \caption{Grid properties for the 4 grids used in the convergence tests, including the number of degrees of freedom (DOFs) for velocity $\MM{u}$, depth $D$ and potential vorticity $q$, along with the timestep used for the solid body rotation test (2) and flow over a mountain test (5).}
  \label{table:grid params}
\end{table}

Where shown, normalised errors in a field $q$ are computed as in \citet{williamson1992standard} as

\begin{align}
  L_2(q) &= \frac{\left(\int_\Omega (q - q_T)^2\diff V\right)^\frac{1}{2}}{\left(\int_\Omega q_T^2\diff V\right)^\frac{1}{2}},\\
  L_\infty(q) &= \frac{\max|q - q_T|}{\max|q_T|},
\end{align}
where $q_T$ is the true solution, specified either analytically (as in
the solid body rotation test) or from a high resolution reference
solution (as in the flow over a mountain test).

\subsection{Solid body rotation (Williamson, test case 2)}

\begin{figure}
  \centering
  \includegraphics[width=0.5\linewidth]{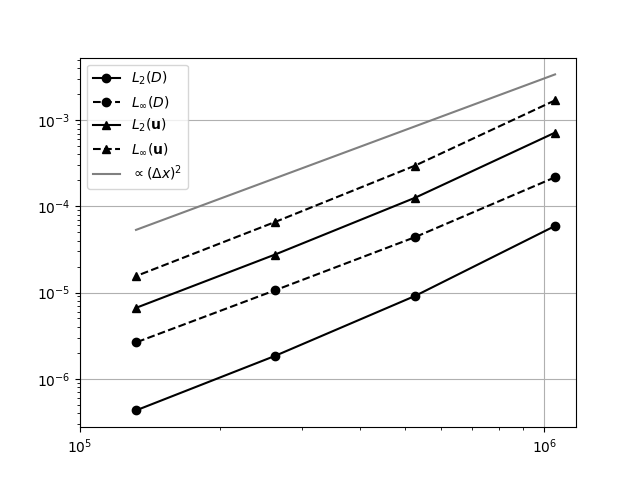}
  \caption{Solid body rotation test case: normalised depth, $D$, and velocity, $\MM{u}$, errors at day {1}5 versus average mesh size $\Delta x$.}
  \label{fig:W2 convergence}
\end{figure}

\begin{figure}
  \centering
  \includegraphics[width=0.7\linewidth]{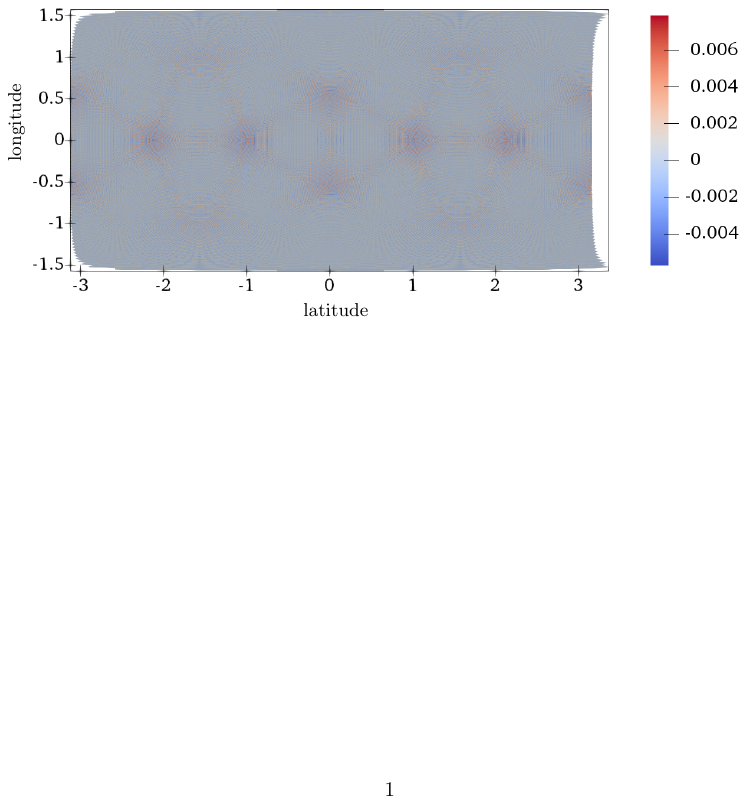}
  \caption{Solid body rotation test case: depth error (in metres) at day 15. The range is min:{\num{-5.8e-3}m}, max: {\num{7.9e-3}m}.}
  \label{fig:W2 D error latlon}
\end{figure}

This test case is initialised with depth and velocity fields that are
in geostrophic balance:

\begin{align}
  D &= D_0 - \left(R\Omega u_0 + \frac{u_0^2}{2}\right)\frac{z^2}{gR^2}, \label{eq:zonal flow D}\\
  u &= \frac{u_0}{R}(-y, x, 0),\label{eq:zonal flow u}
\end{align}

\noindent where $R=6.37122\times 10^6\text{m}$ is the radius of the
Earth, $\Omega=7.292\times 10^{-5}\text{s}^{-1}$ is the rotation rate
of the Earth, $(x, y, z)$ are the $3D$ Cartesian coordinates,
$g=9.80616 \text{ms}^{-2}$ is the gravitational acceleration,
$D_0=2.94\times 10^4/g\approx 2998\text{m}$ and $u_0=2\pi R/(12 \text{
  days})\approx 38.6\text{ms}^{-1}$. Since there is no topography or
other forcing, the flow should remain in this steady, balanced
state. As the flow is steady we have an analytic solution which allows
us to compute errors and hence an order of convergence for our
method. We present these results in table \ref{table:solid body
  errors} and Figure \ref{fig:W2 convergence} reveals that our method
is converging at second order. The depth error at day 15, shown in
Figure \ref{fig:W2 D error latlon}, is large scale and shows some
evidence of grid imprinting.

\begin{table}
  \centering
  \begin{tabular}{rSSSS}
    {\text{cells}} & {$L_2(D)$} & {$L_\infty(D)$} & {$L_2(\MM{u})$} & {$L_\infty(\MM{u})$} \\
    1280 & \num{5.929e-5} & \num{2.177e-4} & \num{7.180e-4} & \num{1.701e-3} \\
    5120 & \num{9.154e-6} & \num{4.405e-5} & \num{1.261e-4} & \num{2.978e-4} \\
    20480 & \num{1.840e-6} & \num{1.062e-5} & \num{2.761e-5} & \num{6.588e-5} \\
    81920 & \num{4.288e-7} &  \num{2.636e-6} & \num{6.640e-6} & \num{1.551e-5} \\
  \end{tabular}
  \caption{Normalised depth, $D$, and velocity, $\MM{u}$, errors at day {1}5 for the solid body rotation test case.}
  \label{table:solid body errors}
\end{table}

\subsection{Flow over a mountain (Williamson, test case 5)}

\begin{figure}
  \centering
  \includegraphics[width=0.7\linewidth]{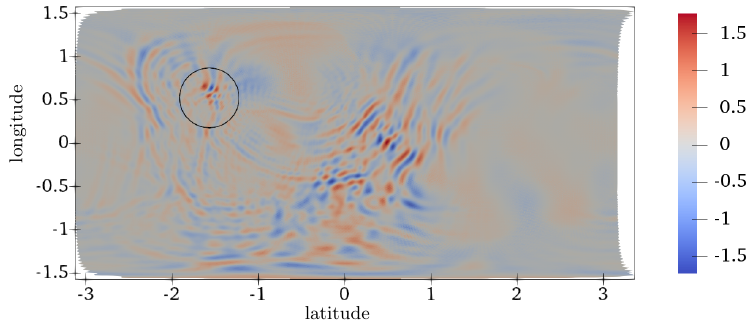}
  \caption{Flow over a mountain test case: depth errors (in metres) at day 15. The range is min:{$-1.74\text{m}$}, max: {$1.78\text{m}$}.}
  \label{fig:W5 D error}
\end{figure}

This test case uses the same zonal flow initial conditions
\eqref{eq:zonal flow D}-\eqref{eq:zonal flow u} as the previous test
but with $D_0=5960\text{m}$ and $u_0=20\text{m}$. An isolated, conical
mountain, given by

\begin{equation}
D_s = D_{s0}(1-r/R_0)
\end{equation}

\noindent with $D_{s0}=2000\text{m}$, $R_0=\pi/9$ and $r^2=\min[R_0^2,
  (\phi-\phi_c)^2+(\lambda-\lambda_c)^2]$ is placed with its centre at
latitude $\phi=\pi/6$ and longitude $\lambda=-\pi/2$. As the zonal
flow interacts with the mountain, it produces waves that travel around
the globe. As there is no analytical solution for this problem, the
model output at 15 days is compared to a high resolution
{(a 2048 by 1024 grid, with a timestep of 45s)} reference
solution generated using a semi-Lagrangian shallow water code provided
by John Thuburn. We plot the depth errors at day 15 in Figure
\ref{fig:W5 D error}. We can see that the error is small scale and is
not dominated by errors due to grid imprinting.

\begin{table}
  \centering
  \begin{tabular}{rSS}
    {\text{cells}} & {$L_2(D)$} & {$L_\infty(D)$} \\
    1280 & \num{1.406e-3} & \num{7.963e-3}  \\
    5120 & \num{6.776e-4} & \num{3.317e-3}  \\
    20480 & \num{2.136e-4} & \num{1.200e-3}  \\
    81920 & \num{4.159e-5} &  \num{2.989e-4}  \\
  \end{tabular}
  \caption{Normalised depth errors at day 15 for the flow over a mountain test case.}
\end{table}

\begin{figure}
  \centering
  \includegraphics[width=0.5\linewidth]{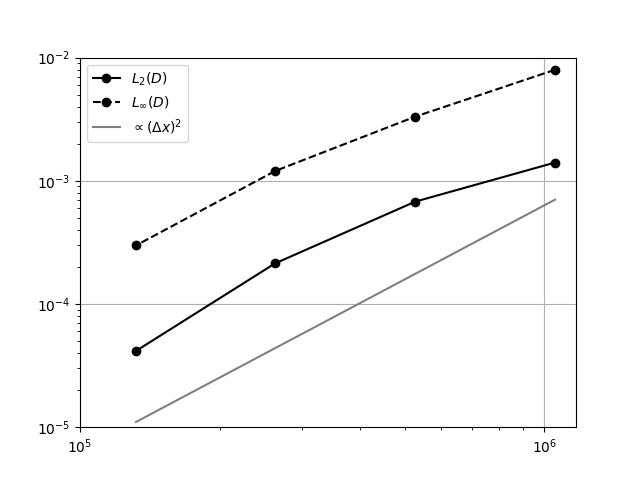}
  \caption{Flow over a mountain test case: normalised depth errors at day 15 versus average mesh size $\Delta x$.}
  \label{fig:W5 convergence}
\end{figure}

\begin{figure}
  \centering
  \includegraphics[width=0.7\linewidth]{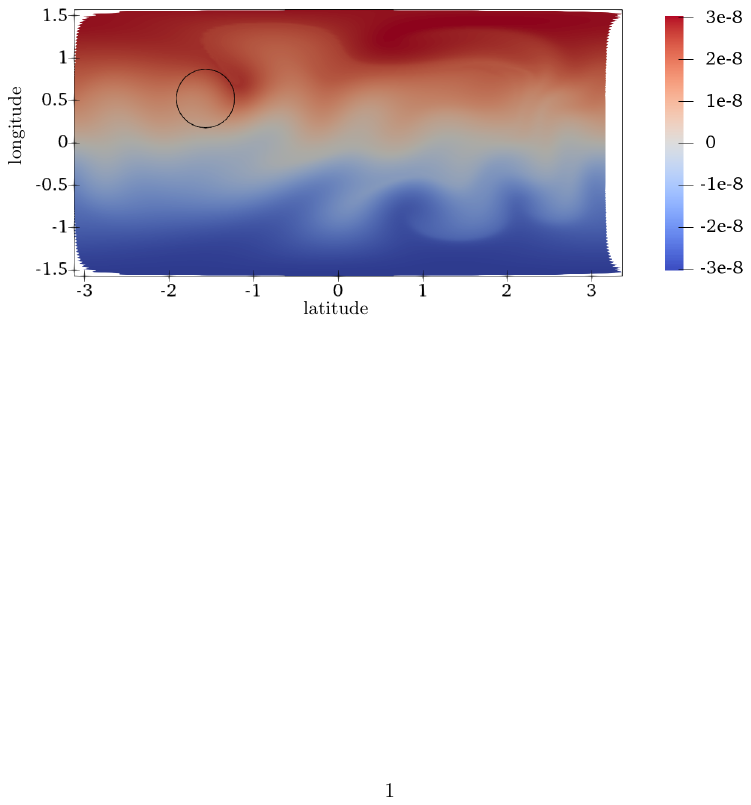}
  \caption{Flow over a mountain test case: potential vorticity at day 15. The range is min:$-3.05\times10^{-8}(\text{ms})^{-1}$, max: $3.05\times10^{-8}(\text{ms})^{-1}$. The circle indicates the position of the mountain.}
  \label{fig:W5 PV day 15}
  \includegraphics[width=0.7\linewidth]{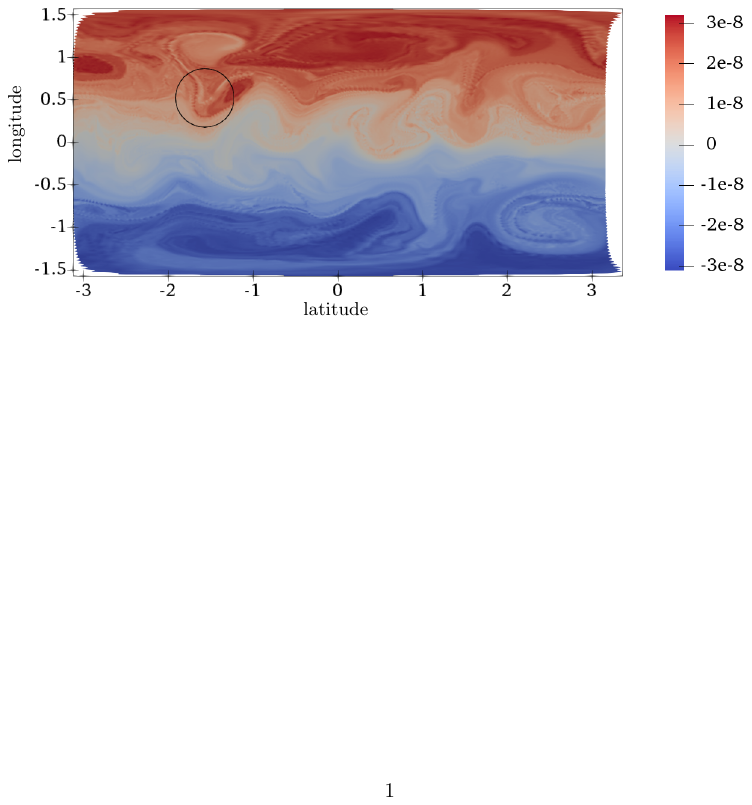}
  \caption{Flow over a mountain test case: potential vorticity at day 50. The range is min:{$-3.1\times10^{-8}(\text{ms})^{-1}$}, max: $3.2\times10^{-8}(\text{ms})^{-1}$. The circle indicates the position of the mountain.}
  \label{fig:W5 PV day 50}
\end{figure}

Up until this point the flow is only weakly nonlinear (see Figure
\ref{fig:W5 PV day 15}) so, as in \citet{thuburn2013mimetic}, we
continued the highest resolution simulation until day 50. By this
time, fine scale structure has been generated as the flow becomes more
nonlinear; the PV field develops sharp gradients and filaments that
stretch out and roll up. These features can be seen in the potential
vorticity field at 50 days, shown in Figure \ref{fig:W5 PV day
  50}.

\subsection{Barotropically unstable jet (Galewsky)}
The details of this test case are specified in
\citet{galewsky2004initial}. The initial condition consists of a
strong midlatitude jet, with an added perturbation, which is
barotropically unstable and evolves to produce vortices and small
scale structure. It has become a particularly useful test for models
on grids that are not aligned with latitude/longitude because these
grids can induce early development of the instability and even lead to
the final solution (at day 6) having the incorrect wavenumber. Again,
there is no analytic solution so results are compared to those in the
literature (see for example \citet{thuburn2013mimetic} and
\citet{weller2013non}). An important feature to reproduce is the
relatively straight path of the jet across approximately quarter of
the globe - this is not seen in models where grid imprinting has
resulted in the generation of instability along the length of the jet
\citep{weller2013non}.

\begin{figure}
  \begin{center}
    \includegraphics[width=0.95\textwidth]{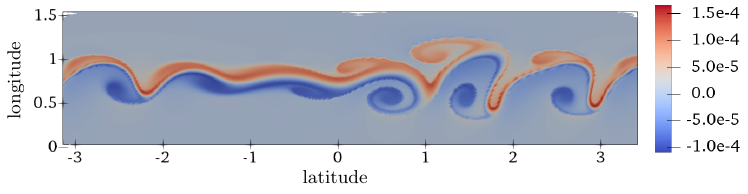}
  \caption{Snapshot of the northern hemisphere relative vorticity
    field at day 6 for the Galewsky jet test. The range is min:
    {$-1.10\times 10^{-4}(\text{ms})^{-1}$}, max:
    {$1.66\times 10^{-4}(\text{ms})^{-1}$}.}
  \label{fig:Galewsky_vorticity}
  \end{center}
\end{figure}

Figure \ref{fig:Galewsky_vorticity} shows the vorticity field,
computed on the highest resolution grid (see table \ref{table:grid
  params}) with timestep of $120\text{s}$, in the Northern hemisphere
after 6 days. The maximum Courant number reached during this
simulation is $0.28$. We note that the jet has the correct wavenumber
and the expected quiescent section. Lower resolution simulations, as
others have found \citep{thuburn2013mimetic}, did not have these
requisite features. \correction{The conservation properties of the
  algorithm are demonstrated in table \ref{table:PV conservation}
  which gives values of the normalised integral of the potential
  vorticity over the sphere,
  \begin{equation}
    \label{eq: Q}
    Q = \frac{\int_\Omega qD \diff V}{\|q_0\|_{L^2}\|D_0\|_{L^2}}
  \end{equation}
  where $q_0$ and $D_0$ are the initial values of the potential
  vorticity and depth respectively, and $\|\cdot\|_{L^2}$ indicates the (un-normalised) $L^2$ norm. This is conserved (and zero) by
  construction (see equation \ref{eq: int qD} and the following
  discussion).}

\begin{table}
  \centering
  \begin{tabular}{cS@{\hskip 1in}rS}
    Time (days) & {$Q$} & & {$Q$} \\
    0 & -5.70e-15 & Min: & -1.80e-14 \\
    2 & -5.22e-15 & Max: & 1.73e-14 \\
    4 & 1.05e-14 & Mean: & 4.23e-16 \\
    6 & -2.02e-16 & Standard deviation: & 4.36e-15
  \end{tabular}
  \caption{Values of $Q$ (see equation \ref{eq: Q}) for the Galewsky test case.}
  \label{table:PV conservation}
\end{table}

\section{Summary and Outlook}
\label{sec: summary}

We have built upon the work of \citet{cotter2012mixed} and
\citet{mcrae2014energy} to produce a semi-implicit compatible finite
element model for the nonlinear rotating shallow water equations on
the sphere. The important features are that we introduce higher-order
upwind advection schemes that maintain PV conservation, and
consistency of PV advection with the mass conservation law. By
applying the model to standard test cases we have demonstrated that
this model has the expected second order convergence rate and can
produce the required features of both large scale balanced flows and
unstable turbulent flows. The developments of this paper inform our
ongoing development of a three dimensional dynamical core on the
sphere.

\paragraph{Acknowledgements} The authors would like to thank John Thuburn for
providing the semi-Lagrangian code that was used for the Williamson 5
convergence test, and the Firedrake project \citep{Rathgeber2016,
  Luporini2016, Homolya2017, Kirby2017, gibson2018domain}, along with PETSc
\citep{petsc-user-ref, petsc-efficient} and other upstream tools
\citep{Chaco95, Dalcin2011} for making the code development in this
paper possible. This research was conducted with support from Natural
Environment Research Council grants NE/I000747/1 and NE/K006789/1.

\bibliographystyle{plainnat}
\bibliography{bibliography}

\begin{thebibliography}{36}
\providecommand{\natexlab}[1]{#1}
\providecommand{\url}[1]{\texttt{#1}}
\expandafter\ifx\csname urlstyle\endcsname\relax
  \providecommand{\doi}[1]{doi: #1}\else
  \providecommand{\doi}{doi: \begingroup \urlstyle{rm}\Url}\fi

\bibitem[Aln{\ae}s et~al.(2014)Aln{\ae}s, Logg, {\O}lgaard, Rognes, and
  Wells]{UFL}
Martin~S. Aln{\ae}s, Anders Logg, Kristian~B. {\O}lgaard, Marie~E. Rognes, and
  Garth~N. Wells.
\newblock Unified form language: A domain-specific language for weak
  formulations of partial differential equations.
\newblock \emph{ACM Transactions on Mathematical Software (TOMS)}, 40\penalty0
  (2):\penalty0 9, 2014.

\bibitem[Arakawa and Lamb(1981)]{arakawa1981potential}
Akio Arakawa and Vivian~R. Lamb.
\newblock A potential enstrophy and energy conserving scheme for the shallow
  water equations.
\newblock \emph{Monthly Weather Review}, 109\penalty0 (1):\penalty0 18--36,
  1981.

\bibitem[Arnold et~al.(2006)Arnold, Falk, and Winther]{arnold2006finite}
DN~Arnold, RS~Falk, and R~Winther.
\newblock Finite element exterior calculus, homological techniques, and
  applications.
\newblock \emph{Acta {N}umerica}, 15\penalty0 (1):\penalty0 1--155, 2006.

\bibitem[Arnold et~al.(2014)Arnold, Boffi, and Bonizzoni]{Arnold14}
DN~Arnold, D~Boffi, and F~Bonizzoni.
\newblock Finite element differential forms on curvilinear cubic meshes and
  their approximation properties.
\newblock \emph{Numer. Math.}, 2014.
\newblock \doi{10.1007/s00211-014-0631-3}.
\newblock arXiv:1204.2595.

\bibitem[Balay et~al.(1997)Balay, Gropp, McInnes, and Smith]{petsc-efficient}
Satish Balay, William~D. Gropp, Lois~Curfman McInnes, and Barry~F. Smith.
\newblock Efficient management of parallelism in object oriented numerical
  software libraries.
\newblock In E.~Arge, A.~M. Bruaset, and H.~P. Langtangen, editors,
  \emph{Modern Software Tools in Scientific Computing}, pages 163--202.
  Birkh{\"{a}}user Press, 1997.

\bibitem[Balay et~al.(2016)Balay, Abhyankar, Adams, Brown, Brune, Buschelman,
  Dalcin, Eijkhout, Gropp, Kaushik, Knepley, McInnes, Rupp, Smith, Zampini,
  Zhang, and Zhang]{petsc-user-ref}
Satish Balay, Shrirang Abhyankar, Mark~F. Adams, Jed Brown, Peter Brune, Kris
  Buschelman, Lisandro Dalcin, Victor Eijkhout, William~D. Gropp, Dinesh
  Kaushik, Matthew~G. Knepley, Lois~Curfman McInnes, Karl Rupp, Barry~F. Smith,
  Stefano Zampini, Hong Zhang, and Hong Zhang.
\newblock {PETS}c users manual.
\newblock Technical Report ANL-95/11 - Revision 3.7, Argonne National
  Laboratory, 2016.

\bibitem[Bochev and Ridzal(2008)]{bochev2008rehabilitation}
PB~Bochev and D~Ridzal.
\newblock Rehabilitation of the lowest-order {R}aviart-{T}homas element on
  quadrilateral grids.
\newblock \emph{SIAM Journal on Numerical Analysis}, 47\penalty0 (1):\penalty0
  487--507, 2008.

\bibitem[Boffi and Gastaldi(2009)]{boffi2009some}
D~Boffi and L~Gastaldi.
\newblock Some remarks on quadrilateral mixed finite elements.
\newblock \emph{Computers \& Structures}, 87\penalty0 (11):\penalty0 751--757,
  2009.

\bibitem[Boffi et~al.(2013)Boffi, Brezzi, and Fortin]{boffi2013mixed}
D~Boffi, F~Brezzi, and M~Fortin.
\newblock \emph{Mixed finite element methods and applications}.
\newblock Springer, 2013.

\bibitem[Brooks and Hughes(1982)]{brooks1982streamline}
Alexander~N Brooks and Thomas~JR Hughes.
\newblock Streamline upwind/{Petrov-Galerkin} formulations for convection
  dominated flows with particular emphasis on the incompressible
  {Navier-Stokes} equations.
\newblock \emph{Computer methods in applied mechanics and engineering},
  32\penalty0 (1):\penalty0 199--259, 1982.

\bibitem[Cotter and Shipton(2012)]{cotter2012mixed}
CJ~Cotter and J~Shipton.
\newblock Mixed finite elements for numerical weather prediction.
\newblock \emph{Journal of Computational Physics}, 231\penalty0 (21):\penalty0
  7076--7091, 2012.

\bibitem[Cotter and Thuburn(2014)]{cotter2014finite}
CJ~Cotter and J~Thuburn.
\newblock A finite element exterior calculus framework for the rotating
  shallow-water equations.
\newblock \emph{Journal of Computational Physics}, 257:\penalty0 1506--1526,
  2014.

\bibitem[Dalcin et~al.(2011)Dalcin, Paz, Kler, and Cosimo]{Dalcin2011}
Lisandro~D. Dalcin, Rodrigo~R. Paz, Pablo~A. Kler, and Alejandro Cosimo.
\newblock Parallel distributed computing using {P}ython.
\newblock \emph{Advances in Water Resources}, 34\penalty0 (9):\penalty0
  1124--1139, 2011.
\newblock \doi{http://dx.doi.org/10.1016/j.advwatres.2011.04.013}.
\newblock New Computational Methods and Software Tools.

\bibitem[Donea(1984)]{donea1984taylor}
Jean Donea.
\newblock A {Taylor--G}alerkin method for convective transport problems.
\newblock \emph{International Journal for Numerical Methods in Engineering},
  20\penalty0 (1):\penalty0 101--119, 1984.

\bibitem[Galewsky et~al.(2004)Galewsky, Scott, and
  Polvani]{galewsky2004initial}
J~Galewsky, RK~Scott, and LM~Polvani.
\newblock An initial-value problem for testing numerical models of the global
  shallow-water equations.
\newblock \emph{Tellus A}, 56\penalty0 (5):\penalty0 429--440, 2004.

\bibitem[Gibson et~al.(2018)Gibson, Mitchell, Ham, and
  Cotter]{gibson2018domain}
Thomas~H Gibson, Lawrence Mitchell, David~A Ham, and Colin~J Cotter.
\newblock A domain-specific language for the hybridization and static
  condensation of finite element methods.
\newblock \emph{arXiv preprint arXiv:1802.00303}, 2018.

\bibitem[Gopalakrishnan and Tan(2009)]{gopalakrishnan2009convergent}
Jayadeep Gopalakrishnan and Shuguang Tan.
\newblock A convergent multigrid cycle for the hybridized mixed method.
\newblock \emph{Numerical Linear Algebra with Applications}, 16\penalty0
  (9):\penalty0 689--714, 2009.

\bibitem[Gottlieb et~al.(2001)Gottlieb, Shu, and Tadmor]{gottlieb2001strong}
Sigal Gottlieb, Chi-Wang Shu, and Eitan Tadmor.
\newblock Strong stability-preserving high-order time discretization methods.
\newblock \emph{SIAM review}, 43\penalty0 (1):\penalty0 89--112, 2001.

\bibitem[Hendrickson and Leland(1995)]{Chaco95}
Bruce Hendrickson and Robert Leland.
\newblock A multilevel algorithm for partitioning graphs.
\newblock In \emph{Supercomputing '95: Proceedings of the 1995 ACM/IEEE
  Conference on Supercomputing (CDROM)}, page~28, New York, 1995. ACM Press.
\newblock ISBN 0-89791-816-9.
\newblock \doi{http://doi.acm.org/10.1145/224170.224228}.

\bibitem[Holst and Stern(2012)]{Holst12}
M~Holst and A~Stern.
\newblock Geometric variational crimes: Hilbert complexes, finite element
  exterior calculus, and problems on hypersurfaces.
\newblock \emph{Foundations of Computational Mathematics}, 12\penalty0
  (3):\penalty0 263--293, 2012.

\bibitem[Homolya et~al.(2017)Homolya, Mitchell, Luporini, and Ham]{Homolya2017}
Mikl{\'o}s Homolya, Lawrence Mitchell, Fabio Luporini, and David~A. Ham.
\newblock {TSFC: a structure-preserving form compiler}, 2017.
\newblock URL \url{http://arxiv.org/abs/1705.003667}.

\bibitem[Kirby and Mitchell(2017)]{Kirby2017}
Robert~C. Kirby and Lawrence Mitchell.
\newblock {Solver composition across the PDE/linear algebra barrier}.
\newblock \emph{SIAM Journal on Scientific Computing}, 2017.
\newblock URL \url{http://arxiv.org/abs/1706.01346}.
\newblock to appear.

\bibitem[Luporini et~al.(2016)Luporini, Ham, and Kelly]{Luporini2016}
Fabio Luporini, David~A. Ham, and Paul H.~J. Kelly.
\newblock An algorithm for the optimization of finite element integration
  loops.
\newblock \emph{Submitted to ACM TOMS}, 2016.
\newblock URL \url{http://arxiv.org/abs/1604.05872}.

\bibitem[McRae and Cotter(2014)]{mcrae2014energy}
Andrew~TT McRae and Colin~J Cotter.
\newblock Energy-and enstrophy-conserving schemes for the shallow-water
  equations, based on mimetic finite elements.
\newblock \emph{Quarterly Journal of the Royal Meteorological Society},
  140\penalty0 (684):\penalty0 2223--2234, 2014.

\bibitem[Natale et~al.(2016)Natale, Shipton, and Cotter]{natale2016compat}
A.~Natale, J.~Shipton, and C.~J. Cotter.
\newblock Compatible finite element spaces for geophysical fluid dynamics.
\newblock \emph{Dyn. Stat. Climate Sys.}, 2016.

\bibitem[Putman and Lin(2007)]{putman2007finite}
WM~Putman and SJ~Lin.
\newblock Finite-volume transport on various cubed-sphere grids.
\newblock \emph{Journal of Computational Physics}, 227\penalty0 (1):\penalty0
  55--78, 2007.

\bibitem[Rathgeber et~al.(2016)Rathgeber, Ham, Mitchell, Lange, Luporini,
  McRae, Bercea, Markall, and Kelly]{Rathgeber2016}
Florian Rathgeber, David~A. Ham, Lawrence Mitchell, Michael Lange, Fabio
  Luporini, Andrew T.~T. McRae, Gheorghe-Teodor Bercea, Graham~R. Markall, and
  Paul H.~J. Kelly.
\newblock Firedrake: automating the finite element method by composing
  abstractions.
\newblock \emph{ACM Trans. Math. Softw.}, 43\penalty0 (3):\penalty0
  24:1--24:27, 2016.
\newblock ISSN 0098-3500.
\newblock \doi{10.1145/2998441}.
\newblock URL \url{http://arxiv.org/abs/1501.01809}.

\bibitem[Ringler et~al.(2010)Ringler, Thuburn, Klemp, and
  Skamarock]{ringler2010unified}
T.~D. Ringler, J.~Thuburn, J.~B. Klemp, and W.~C. Skamarock.
\newblock A unified approach to energy conservation and potential vorticity
  dynamics for arbitrarily-structured {C}-grids.
\newblock \emph{Journal of Computational Physics}, 229\penalty0 (9):\penalty0
  3065--3090, 2010.
\newblock \doi{10.1016/j.jcp.2009.12.007}.

\bibitem[Rognes et~al.(2013)Rognes, Ham, Cotter, and
  McRae]{rognes2013automating}
ME~Rognes, DA~Ham, CJ~Cotter, and ATT McRae.
\newblock Automating the solution of {PDE}s on the sphere and other manifolds
  in {FEniCS} 1.2.
\newblock \emph{Geoscientific Model Development Discussions}, 6\penalty0
  (3):\penalty0 3557--3614, 2013.

\bibitem[Safjan and Oden(1993)]{safjan1993high}
A~Safjan and JT~Oden.
\newblock High-order {T}aylor-{G}alerkin and adaptive h-p methods for
  second-order hyperbolic systems: Application to elastodynamics.
\newblock \emph{Computer Methods in Applied Mechanics and Engineering},
  103\penalty0 (1):\penalty0 187--230, 1993.

\bibitem[Staniforth and Thuburn(2012)]{staniforth2012horizontal}
Andrew Staniforth and John Thuburn.
\newblock Horizontal grids for global weather and climate prediction models: a
  review.
\newblock \emph{Quarterly Journal of the Royal Meteorological Society},
  138\penalty0 (662):\penalty0 1--26, 2012.

\bibitem[Thuburn and Cotter(2012)]{thuburn2012framework}
J~Thuburn and CJ~Cotter.
\newblock A framework for mimetic discretization of the rotating shallow-water
  equations on arbitrary polygonal grids.
\newblock \emph{SIAM Journal on Scientific Computing}, 34\penalty0
  (3):\penalty0 B203--B225, 2012.

\bibitem[Thuburn et~al.(2013)Thuburn, Cotter, and Dubos]{thuburn2013mimetic}
J~Thuburn, CJ~Cotter, and T~Dubos.
\newblock A mimetic, semi-implicit, forward-in-time, finite volume shallow
  water model: comparison of hexagonal-icosahedral and cubed sphere grids.
\newblock \emph{Geoscientific Model Development Discussions}, 6:\penalty0
  6867--6925, 2013.

\bibitem[Thuburn and Cotter(2015)]{thuburn2015primal}
John Thuburn and Colin~J Cotter.
\newblock A primal--dual mimetic finite element scheme for the rotating shallow
  water equations on polygonal spherical meshes.
\newblock \emph{Journal of Computational Physics}, 290:\penalty0 274--297,
  2015.

\bibitem[Weller(2013)]{weller2013non}
Hilary Weller.
\newblock Non-orthogonal version of the arbitrary polygonal c-grid and a new
  diamond grid.
\newblock \emph{Geoscientific Model Development}, 7:\penalty0 779--797, 2013.

\bibitem[Williamson et~al.(1992)Williamson, Drake, Hack, Jakob, and
  Swarztrauber]{williamson1992standard}
DL~Williamson, JB~Drake, JJ~Hack, R~Jakob, and PN~Swarztrauber.
\newblock A standard test set for numerical approximations to the shallow water
  equations in spherical geometry.
\newblock \emph{Journal of Computational Physics}, 102\penalty0 (1):\penalty0
  211--224, 1992.

\end{thebibliography}

\end{document}